\documentclass[notitlepage]{amsart}

\usepackage{amsfonts}
\usepackage[pdftex]{graphicx}
\usepackage[latin1]{inputenc}
\usepackage[all]{xy} 
\usepackage{latexsym,amssymb,amscd}
\usepackage{amsmath}
\usepackage{color}
\newtheorem{pro}{Proposition}[section]
\newtheorem{teo}[pro]{Theorem}
\newtheorem{defi}[pro]{Definition}
\newtheorem{lem}[pro]{Lemma}
\newtheorem{cor}[pro]{Corollary}
\newtheorem{rk}[pro]{Remark}

\newtheorem{ex}[pro]{Example}

\newcommand{\I}{\mathcal{I}}
\newcommand{\Pro}{\mathcal{P}}

\newcommand{\C}{\mathcal{C}}

\newcommand{\F}{\mathfrak{F}}
\newcommand{\A}{\mathcal{A}}

\newcommand{\X}{\mathcal{X}}

\newcommand{\modu}{\mathrm{mod}}

\newcommand{\proj}{\mathrm{proj}}

\newcommand{\Ext}{\mathrm{Ext}}
\newcommand{\Hom}{\mathrm{Hom}}
\newcommand{\End}{\mathrm{End}}
\newcommand{\Ker}{\mathrm{Ker}}
\newcommand{\Coker}{\mathrm{Coker}}
\newcommand{\Ima}{\mathrm{Im}}
\newcommand{\pd}{\mathrm{pd}}
\newcommand{\gld}{\mathrm{gl.dim}}

\newcommand{\Tr}{\mathrm{Tr}}
\newcommand{\rad}{\mathrm{rad}}
\newcommand{\topo}{\mathrm{top}}
\newcommand{\soc}{\mathrm{soc}}

\newcommand{\add}{\mathrm{add}}
\newcommand{\ann}{\mathrm{ann}}
\newcommand{\pres}{\mathrm{pres}}
\newcommand{\wide}{\mathrm{wide}}

\newenvironment{dem}{\noindent\bf Proof. \rm }{$\ \Box$}
\usepackage{latexsym,amssymb,amscd}
\usepackage{amsmath}

\begin{document}
\title[wide subcategories]{Wide subcategories of finitely generated $\Lambda$-modules.}
\author{E.N.Marcos, O. Mendoza, C. S\'aenz, V.Santiago.}
\thanks{
The authors thanks the Project PAPIIT-Universidad Nacional Aut\'onoma de M\'exico IN102914.\\The first named author  was partially supported by Conselho Nacional de Desenvolvimento Cient\'{\i}fico
Tecnol\'{o}gico (CNPq, Brazil)  and by the Projeto Fapesp 2014/09310-5.}
\date{}
\begin{abstract}  We explore some properties of wide subcategories of the category 
$\modu\,(\Lambda)$ of finitely generated left $\Lambda$-modules, for some artin algebra $\Lambda.$
In particular we look at wide finitely generated subcategories  and give a connection with the class of standard modules and standardly  stratified algebras. Furthermore, for a wide class $\X$ in $\modu\,(\Lambda),$ we give necessary and sufficient conditions to see that 
$\X=\pres\,(P),$ for some projective $\Lambda$-module $P;$  and finally, a connection with ring epimorphisms is given. 
\end{abstract}
\maketitle
\section*{Introduction.} 
Wide categories have been extensively studied in the case of modules over commutative rings \cite{H,K,T}. One important feature 
of them is their connection with a classification theorem of the thick subcategories of the derived category of perfect complexes over 
a commutative noetherian ring in terms of the ring spectrum.  For non-commutative rings, as far as we know, there are some 
interesting results in  \cite{IT, IPT} for the case of wide subcategories of finitely generated left $A$-modules, where $A$ is 
a finite dimensional, basic and hereditary algebra over an algebraically closed field $k.$ On the other hand, in \cite{MS} there are stated some nice correspondences between torsion classes and wide classes in   $\modu\,(A),$ for a finite dimensional $k$-algebra. 
 Furthermore, in the case that $A$ is of finite representation type, explicit bijections are given in \cite{MS}  between different classes and wide 
 classes in   $\modu\,(A).$ It seems interesting to us to look at what happens if we consider the more general case of an artin algebra 
 $\Lambda.$ This paper goes in this direction, we study some elementary properties of wide subcategories of $\modu\,(\Lambda),$ 
 for an artin algebra $\Lambda,$ where $\modu\,(\Lambda)$ is the full subcategory of  finitely generated left $\Lambda$-modules. 
\

Let $\A$ be an abelian category. A full subcategory $\C\subseteq\A$ is wide if $\C$ is closed under extensions, kernels and cokernels  of  morphisms in $\C.$ Thus, a wide subcategory is  nothing more than an abelian subcategory which is closed under extensions. Let $\X$ be a class of objects in $\modu\,(\Lambda).$ We denote by $\wide\,(\X)$ the smallest wide subcategory of $\modu\,(\Lambda)$ containing $\X.$ In this paper we study some special types of wide subcategories of $\modu\,(\Lambda).$
\

The paper is organized as follows.  In section 2, we give some basic notions and an equivalent description for a subcategory to be wide.
\


Let $M$ be an indecomposable $\Lambda$-module. In section 3, we prove that 
$\add\,(M)=\wide\,(M)$ if and only if $\Ext^1_\Lambda(M,M)=0$ and $\End_\Lambda(M)$ is a division ring. We also give examples where $\add\,(M)\subsetneq \wide\,(M),$ but still $\wide\,(M)=\add\,(N)$ for some other $\Lambda$-module $N.$ 
\

In section 4, we deal with the set $\Delta$ of standard $\Lambda$-modules and the class $\F(\Delta)$ of good $\Lambda$-modules.  Here, we prove that $\wide\,(\Delta)=\modu\,(\Lambda),$ for any artin $R$-algebra $\Lambda.$  Some consequences, of the above result, are also given in this section. In particular, for a quotient path algebra $\Lambda=kQ/I,$ we prove that the quiver $Q$ is directed, in the sense of S. Koenig \cite{K}, if and only if $\F(\Delta)$ is wide.  
\

In section 5, we deal with the category $\pres\,(P)$ of $\Lambda$-modules having a presentation in 
$\add\,(P),$ for some projective $\Lambda$-module $P.$ It is well known that $\pres\,(P)$ is an abelian category, but in general, it is not true (an example is given) that $\pres\,(P)$ is an abelian subcategory of $\modu\,(\Lambda).$ In this section, we characterize wide subcategories $\X$ of the form  $\X=\pres\,(P)$ for some $P\in\proj\,(\Lambda).$ We prove that, for  a wide subcategory $\X$ of 
$\modu\,(\Lambda),$  the following statements are equivalent
\begin{itemize}
\item[(a)] $\X=\pres\,(P)$ for some $P\in\proj\,(\Lambda);$
\item[(b)] $\X$ has enough projectives and $\proj\,(\X)\subseteq \proj\,(\Lambda);$
\item[(c)] $\X$ is closed under projective covers, that is $\mathbb{P}_0(\X)\subseteq\X.$
\end{itemize}
If one of the above equivalent conditions hold, then 
$$\add\,(P)=\proj\,(\X)=\add\,(\mathbb{P}_0(\X)).$$ 
In particular, $P$ is uniquely determined  up to additive closures. Furthermore,  the functor $\Hom_\Lambda(P,-):\X\to 
\modu\,(\End_\Lambda(P)^{op})$ is an equivalence of categories.
\

Finally, in the case of finite dimensional $k$-algebras and assuming that $\X$ is a  functorially finite wide subcategory  of $\modu\,(A),$ we have that the following statements are equivalent
\begin{itemize}
\item[(a)] There is a ring epimorphism $A\to B$ such that $B\in\proj\,(A)$ and $\X=\pres\,({}_AB).$
\item[(b)] $\X$ has enough projectives and $\proj\,(\X)\subseteq \proj\,(A).$
\end{itemize}
If one of the above equivalent conditions hold,  there is a basic  $A$-module $P$ such that $\add\,(P)=\proj\,(\X)$ and $B$ is 
Morita equivalent to $\End_A(P)^{op}.$

\section{Preliminaries.}

Throughout this paper, the term algebra means {\it artin algebra} over a commutative artin 
ring $R.$ For an algebra $\Lambda,$ the category of finitely generated left $\Lambda$-modules will be
denoted by $\modu\,(\Lambda).$ Unless otherwise specified, we will work  with finitely generated $\Lambda$-modules,
full subcategories and non-empty classes. We denote by $\proj\,(\Lambda)$ the full subcategory of 
$\modu\,(\Lambda)$ whose objects are the projective $\Lambda$-modules. 

Let us denote by $S_1, S_2,\cdots ,S_n$  a complete list of non-isomorphic simple $\Lambda$-modules and let $\leq$ denotes this fixed natural ordering of the simple modules. Let $P_i$ be the projective cover of the simple module $S_i,$ for each $i=1,2,\cdots, n.$ We define for each $i,$ the standard module 
$\Delta(i),$ which is the quotient $\Lambda$-module $\Delta(i)=P_i/U_i$ where $U_i$ is the trace of $\oplus_{j>i}P_j$ in $P_i.$ That is, $U_i=\sum_{f:\oplus_{j>i}P_j\rightarrow P_i}\mathrm{Im} f$.
Let $\Delta=\{\Delta(1),\Delta(2),\cdots,\Delta(n)\}$ be the set of standard modules. It is well known, that the standard module $\Delta(i)$ is the maximal quotient of $P_i$ which composition factors are only the simple $\Lambda$-modules $S_j$ with $j\leq i.$ 
\

For a given class $\Theta$ of $\Lambda$-modules,  we denote by $\F(\Theta)$ the full subcategory of $\modu\,(\Lambda),$ whose objects are  the $\Lambda$-modules $M$ which have a $\Theta$-filtration. That is, $M\in\F(\Theta)$ if there is a finite chain
$$0=M_0\leq M_1\leq M_2<\cdots \leq M_t=M$$

\noindent of submodules of $M$ such that each quotient $M_i/M_{i-1}$ is isomorphic to a module in $\Theta.$ In case of the class $\Theta:=\Delta,$ the modules in $\F(\Delta)$ are called good modules.  

\noindent If ${}_\Lambda\Lambda\in\F(\Delta)$ then the pair $(\Lambda,\leq)$ is said to be a (left) standardly stratified algebra. A standardly stratified algebra is called quasi-hereditary if the endomorphism ring of each standard module is a division ring. Quasi-hereditary algebras were introduced by L.L. Scott in \cite{S}.

Associated with the category  $\F(\Delta),$  there is the class of the $\Delta$-injective $\Lambda$-modules
$$\I(\Delta):=\{Y\in\modu\,(\Lambda)\,:\, \Ext_\Lambda^1(-,Y)|_{\F(\Delta)}=0\} .$$
In the theory of standardly stratified algebras, the subcategory $\F(\Delta)\cap \I(\Delta)$ is of 
special interest.
In fact, this subcategory was characterized by C. M. Ringel in \cite {R} and by I.Agoston, D. Happel, E. Lukacs and L. Unger in \cite {AHLU}, for quasi-hereditary algebras and for standardly stratified algebras, respectively. Their 
characterization is the following theorem.

\begin{teo}[\cite{R, AHLU}]
Let $(\Lambda,\leq)$ be a  standardly stratified algebra. Then, there is a (generalized) tilting basic module $T$  such that $\add (T)=\F(\Delta)\cap \I(\Delta).$
\end{teo} 

\noindent Finally, we recall that for any $M\in\modu\,(\Lambda),$ the class  $\add (M)$ consists of all $\Lambda$-modules which are direct summands of finite coproducts of  copies of $M.$

\vspace{.3cm}

\section{Wide subcategories}

We start this section introducing some basic definitions, which will be used in the sequel. 
\

Let $\A$ be an abelian category and 
$\C$ be a class in $\A.$ We say that $\C$ is {\bf closed under extensions} if, for any 
exact sequence $0\to X\to Y\to Z\to 0$ in $\A,$ the fact that $X,Z\in\C$ implies that 
$Y\in\C.$ We say that $\C$ is {\bf thick} if $\C$ is closed under direct summands  and for any 
exact sequence $0\to X\to Y\to Z\to 0$ in $\A,$ the fact that two of the terms of the exact sequence belong to 
$\C$ implies that the third one belongs to $\C.$ We say that $C\in\C$ satisfies the {\bf intersection property}  if for any sub-objects $X$ and $Y$ of $C,$  which are in $\C$, we have that 
$X\cap Y\in\C.$ 

\noindent
A subcategory $\C$ of an abelian category $\A$ is called {\bf admissible}   if it is thick and any $C\in\C$ satisfies the intersection property (see \cite{P}). Note that, the fact that $\C$ is admissible in $\A$ does not mean that the inclusion functor $\C\hookrightarrow\A$ admits a right or a left adjoint functor. Therefore,  the usual definition of an admissible subcategory of a triangulated category is not related to the above one.

Following M. Hovey, in \cite{H}, we say that 
the subcategory $\C\subseteq\A$ is {\bf wide} if it is closed under extensions, kernels and cokernels of morphisms in $\C.$ Wide categories have been extensively studied in the case of modules over commutative rings \cite{H,K,T}. For the case of finitely generated modules over a hereditary algebra, we recommend the reader to see \cite{IT,IPT}.
An easy observation is that any wide category $\C$ is closed under isomorphisms and furthermore it is an additive subcategory of $\A.$ 
\
 
Given  a class $\X\subseteq \A$, we denote by $\wide\,(\X)$  the smallest wide subcategory of $\A$ 
containing $\X$. Note that the class $\F(\X)$ of all the $\X$-filtered objects is a subclass of $\wide\,(\X),$ since $\F(\X)$ is the smallest subclass in $\A,$ which is closed under extensions and contains $\X.$
\

Let $\C$ be an abelian category, which is an additive subcategory of an abelian category $\A.$ We recall that $\C$ is an abelian subcategory of $\A$ if the inclusion functor $i:\C\rightarrow \A$ is exact. The following is a nice criteria to check that an additive subcategory is an abelian subcategory.

\begin{rk}\label{criterioA} Let $\C$ be an additive subcategory of an abelian category $\A.$ Then, $\C$ is an abelian subcategory of $\A$ if and only if $\C$ is closed under kernels and cokernels (in $\A$)
 of morphisms in $\C.$
\end{rk}

\begin{cor} \label{CohAb}
Let $\C$ be a wide subcategory of an abelian category $\A.$ Then, $\C$ is abelian and the inclusion functor 
$i:\C\rightarrow \A$ is full, faithful and exact. Moreover, 
$\Ext^1_\C(X,Y)=\Ext^1_\A(X,Y)$ for any $X,Y\in\C.$   
\end{cor}
\begin{dem} It follows from the definition of wide category and Remark \ref{criterioA}.
\end{dem}

\begin{pro}\label{R1}
Let $\C$ be a subcategory of an abelian category $\A$. Then, 
$\C$ is admissible if and only if $\C$ is wide.  
\end{pro}
\begin{dem}
$(\Rightarrow)$ Assume that $\C$ is admissible.
Let us show that $\C$ is closed under kernels and cokernels.

\noindent Let $f:A\rightarrow B$ be a morphism in $\C$. Consider the sub-objects $M$ and $N$ of $A\oplus B\in \C,$ where
 $$M:=\{(a,f(a))\,:\, a\in A\}\,\,\mathrm{and}\,\, N:=\{(a,0)\,:\,a\in A\}.$$ 
 Since $g:A\rightarrow M,$ given by $a\mapsto (a,f(a)),$ is an 
 isomorphism and $N\simeq A,$ we get that $M$ and $N\in \C$. Therefore, using that $A\oplus B\in\C,$ it follows that $\Ker\,(f)\simeq M\cap N\in \C$. Thus, the exact sequence  
$$0\rightarrow \Ker\,(f)\rightarrow A\rightarrow \Ima\,(f)\rightarrow 0$$ gives us that 
$\Ima\,(f)\in\C.$ Finally, by considering the exact sequence  
$$0\rightarrow \Ima\,(f)\rightarrow A\rightarrow \Coker\,(f)\rightarrow 0,$$ we get that 
$\Coker\,(f)\in \C.$
\

$(\Leftarrow)$ Suppose that $\C$ is wide. We assert that $\C$ is closed under direct summands. Indeed,  let $C\in\C$  and consider a decomposition $C=A\oplus B$ in 
$\modu(\Lambda)$ and let $\alpha:=\left(\begin{matrix}
1_A & 0\\0&0 \end{matrix}\right):C\rightarrow C.$ Thus  $B\simeq\Ker\,(\alpha)\in\C,$ proving that $\C$ is closed under direct summands.

Finally, let $M_1$ and $M_2$ be sub-objects of $M$ with $M_1,M_2$ and $M\in \C$. Consider the morphism 
$$f:M_1\oplus M_2\rightarrow M,\quad (m_1,m_2)\mapsto m_1+m_2.$$ 
Then, since $\C$ is closed under kernels of morphisms in $\C,$ we get 
$$M_1\cap M_2\simeq \{(m,-m) \,:\, m\in M_1\cap M_2\}=\Ker\,(f)\in\C;$$ proving that 
$\C$ is admissible.
\end{dem}

\begin{cor}\label{R3}
Let $\C$ be a wide subcategory of $\modu (\Lambda).$ If  ${}_\Lambda\Lambda\in \C$ then $\C=\modu (\Lambda).$
\end{cor}

\begin{dem}
For every $M\in\modu (\Lambda)$ there exists an exact sequence 
$${}_\Lambda\Lambda^r\stackrel{f}{\rightarrow}{}_\Lambda\Lambda^s\rightarrow M\rightarrow 0,$$
and so we have that $M\simeq\Coker f\in\C.$
\end{dem}

\section{Wide finitely generated subcategories }

Let $\Lambda$ be an artin $R$-algebra. For any $M\in\modu\,(\Lambda),$  we sometimes use $(M,-)$ to denote the functor 
$\Hom_\Lambda(M,-).$
\

 Consider  a wide subcategory $\X$ of $\modu\,(\Lambda).$ 
We say that $\X$ is {\bf  wide finitely generated} if there is some $X\in\X$ such that $\X=\wide\,(X).$ We 
also say that  $\X$ is {\bf wide additively generated} if there is some $X\in\X$ such that $\X=\add\,(X).$ 
\

In what follows, for some indecomposable $M\in\modu\,(\Lambda),$
we  give necessary and sufficient conditions
 for the category $\add\,(M)$ to be wide. For doing so, we need the following result.

\begin{pro}\label{R2}
Let $M\in \modu(\Lambda)$ be such that $\End_\Lambda(M)$ is a semi-simple ring. Then, the category  $\add\,(M)$ is closed under kernels, images and cokernels of morphisms in $\add\,(M).$ In particular, 
$\add\,(M)$ is an abelian subcategory of $\modu\,(\Lambda)$ and $(M,-):\add\,(M)\to \modu\,(\Gamma)$ is an equivalence of categories, where $\Gamma:=\End_\Lambda(M)^{op}.$
\end{pro}
\begin{dem} 
It is well known (see \cite{ARS, ASS})  that the evaluation  functor $(M,-):\modu\,(\Lambda)\to \modu\,(\Gamma)$ 
induces an equivalence of $R$-categories between $\add\,(M)$ and $\proj\,(\Gamma).$
\

By hypothesis we know that $\Gamma$ is a semi-simple ring. Thus, we have that 
$(M,-):\add\,(M)\to \modu\,(\Gamma)$ is an equivalence of categories.
\
 
Firstly we assert that $\add\,(M)$ has the following property: For all $X,E\in\add\,(M)$ any monomorphism (epimorphism) $\alpha:X\to E$ is a split-mono (split-epi). In particular, the class $\add\,(M)$ is closed under cokernels (kernels) of monomorphisms (epimorphisms). Indeed, let 
$\alpha: X\to E$ be a monomorphism with $X,E\in\add\,(M).$ Therefore $(M,\alpha):(M,X)\to (M,E)$ is a monomorphism in $\modu\,(\Gamma).$ Thus 
 $(M,\alpha)$ is a split-mono, since $\Gamma$ is a semi-simple ring. But now, using the fact that $(M,-):\add\,(M)\to \modu\,(\Gamma)$ is full and faithful, we get that $\alpha$ is a split-mono.    
\

Let $f:A\rightarrow B$ be a morphism in $\add\,(M).$ We prove that $\Ima\,(f)$ and $\Coker\,(f)$ belong to $\add\,(M).$ Note that, from the exact sequence $0\rightarrow \Ima\,(f)\rightarrow B\rightarrow \Coker\,(f)\rightarrow 0$ and the fact that $\add\,(M)$ is closed under cokernels of monomorphisms, we reduce the problem to prove only that $\Ima\,(f)\in\add\,(M).$  Indeed, consider 
the morphism $(M,f):(M,A)\to (M,B)$ of $\Gamma$-modules and its factorization trough its
image $\Ima\,(M,f).$ Since $(M,-):\add\,(M)\to \modu\,(\Gamma)$ is an equivalence of categories and 
$\Gamma$ is a semi-simple ring, we have that $\Ima\,(M,f)=(M,M')$ and $(M,f):(M,A)\to (M,B)$ is the composition of a split-epi $(M,\alpha):(M,A)\to (M,M')$ and a split-mono $(M,\beta):(M,M')\to (M,B).$ In particular, we get that $f=\beta\alpha,$ where $\alpha$ is a split-epi and $\beta$ is 
a split-mono. Thus $\Ima\,(f)\simeq M'\in\add\,(M).$ Finally, to prove that $\Ker\,(f)\in\add\,(M),$ we use the exact sequence $0\to\Ker\,(f)\to A\to \Ima\,(f)\to 0.$
\end{dem}
\vspace{.2cm}

Given an artin algebra $\Lambda$ and $M\in\modu\,(\Lambda)$, we denote by $\ell_\Lambda(M)$  the  length of a Jordan-H\"older composition series of  $M.$ 

\begin{teo}\label{CR3}
Let $M$ be an indecomposable $\Lambda$-module. Then, 
$\add\,(M)=\wide\,(M)$ if and only if $\Ext^1_\Lambda(M,M)=0$ and $\End_\Lambda(M)$ is a division ring. 
\end{teo}

\begin{dem} $(\Rightarrow)$  Assume that $\add (M)=\wide\,(M).$ Let us see that  $\Ext^1_\Lambda(M,M)=0.$ Indeed, consider an 
exact sequence 
$$\eta:\quad  0\rightarrow M\stackrel{\alpha}{\rightarrow}L\stackrel{\beta}\rightarrow M\rightarrow 0$$
 in $\modu (\Lambda).$
Since $\add (M)$ is closed under extensions and $M$ is indecomposable, we get that $L\simeq M^t.$ Thus  $t\,\ell_\Lambda(M)=\ell_\Lambda(L)=2\,\ell_\Lambda(M)$ and so $L\simeq M^2.$ 

Consider the evaluation functor $(M,-):\modu\,(\Lambda)\to \modu\,(\Gamma),$ where $\Gamma:=\End_\Lambda(M)^{op}.$  By applying the  functor $(M,-)$ to  $\eta,$ we get the exact sequence in $\modu (\Gamma)$
$$(M,\eta):\quad  0\rightarrow (M,M)\xrightarrow{(M,\alpha)}(M,L)\xrightarrow{(M,\beta)}(M, M).$$

We claim that $\Ima\, (M,\beta)=(M,M).$ Indeed, suppose that this is not the case, and so  $\ell_\Gamma\Ima\, (M,\beta)<\ell_\Gamma(M,M).$ 
Furthermore,  the exact sequence of $\Gamma$-modules
$$0\rightarrow(M,M)\to (M,L)\to \Ima\,(M, \beta)\rightarrow 0,$$ 
gives us the following
$$\ell_\Gamma (M,L)=\ell_\Gamma(M,M)+\ell_\Gamma(\Ima\, (M,\beta))<2\,\ell_\Gamma (M,M).$$

\noindent Using now that $L\simeq M^2,$ we get  $\ell_\Gamma (M,L)=2\,\ell_\Gamma(M,M);$ which is a contradiction, 
proving that  $\Ima\, (M,\beta)=(M,M).$ 

\noindent The fact that  $\Ima\, (M,\beta)=(M,M)$ implies that the exact sequence $\eta$ splits and hence $\Ext^1_\Lambda(M,M)=0$
\

We prove now that $\End_\Lambda(M)$ is a division ring. Suppose  that this is not the case, and so the radical of the ring  $\End_\Lambda(M)$ is non-zero. Then,  there exists a non-zero nilpotent endomorphism $f:M\to M.$ Note that $\Ker\,(f)\neq 0,$ since $M$ is non-zero.
\

Since $\add\,(M)$ is wide, we 
get that $\Ker\,(f)\in\add\,(M).$ By using Krull-Schmidt theorem and the fact that $M$ is indecomposable, we conclude that 
$\Ker\,(f)\simeq M^t,$ for some $t\geq 1.$ Therefore $\Ker\,(f)= M$ and then $f=0,$ which is a contradiction; proving that  $\End_\Lambda(M)$ is a division ring
\

$(\Leftarrow)$ It follows from Proposition \ref{R2} and using the fact that $\Ext_\Lambda^1(M,M)=0$ implies that $\add (M)$ is closed under extensions.
\end{dem}

\begin{rk} Let $\Lambda$ be a finitely dimensional $k$-algebra, over an algebraically closed field 
$k.$ 
\begin{itemize}
\item[(1)] Any indecomposable $\Lambda$-module $M,$ which is post-projective or pre-injective, satisfies that 
$\End_\Lambda(M)=k$ and $\Ext^{1}_\Lambda(M,M)=0$ \cite[Lemma VIII.2.7]{ASS}. 
\item[(2)] Let $\Lambda$ be a hereditary algebra of finite representation type. Then, any indecomposable $\Lambda$-module $M$ satisfies that 
$\End_\Lambda(M)=k$ and $\Ext^{1}_\Lambda(M,M)=0$ \cite[Corollary VII.5.14]{ASS}
\end{itemize}
Therefore, by Theorem \ref{CR3} each such module produces an example of a wide additively generated  subcategory of 
$\modu\,(\Lambda).$
\end{rk}

\begin{rk}
In general we could have that $\add\,(M)\subsetneq \wide\,(M),$ but still $\wide\,(M)=\add\,(N)$ for some other $\Lambda$-module $N.$  To see this, take an algebra which is of finite representation type and an indecomposable module which has a self-extension. Then it is clear that $\wide\, (M) $ is of the form $\add\, (N)$, since for finite 
representation type algebras $\Lambda$ every full subcategory which is closed by direct summands is of the type $\add(N)$, just take $N$ the direct summand of the indecomposable in $N$. It is also clear that if $M$ has self-extension, then $\wide\, (M)=\add(N)$ properly contains $\add(M).$

Even if we take $M$ indecomposable without self-extensions it can happen that $\wide(M) =\add(N)$ which
properly contains $\add(M),$ for instances of that,  see the following examples. 
\end{rk}

\begin{ex}
Let $n>2$ and $\Lambda=k \tilde{A}_{n-1}/J^{n+1}$ be a symmetric Nakayama algebra  where $J$ denotes the arrow ideal of the path $k$-algebra $\tilde{A}_{n-1}$ \cite[Proposition V.3.8]{ASS}. Note that  $P_1/\rad^n P_1\simeq P_1/\mathrm{soc}\,(P_1).$   Moreover,
we have that  $$\wide\, (P_1)=\add[(P_1/\rad P_1)\oplus (P_1/\rad^n P_1)\oplus P_1\oplus (P_2/\rad^{n-1}P_2)\oplus (P_2/\rad^nP_2)]$$
Indeed, the reader can check that $$\add[(P_1/\rad P_1)\oplus (P_1/\rad^n P_1)\oplus P_1\oplus (P_2/\rad^{n-1}P_2)\oplus (P_2/\rad^nP_2)]$$ is
closed under extensions, kernels and cokernels of morphisms.
\end{ex}

\begin{ex} Let $\Lambda=k\tilde{A}_{2}/I$  where $\tilde{A}_{2}$ is the quiver 
$$
\xymatrix{
						&2 \ar@{->}[ddr]^\beta 	&\\
						&						&\\
1 \ar@{->}[uur]^\alpha 	&						&3 \ar@{->}[ll]^\gamma
}
$$

and $I=\langle \beta\alpha\rangle.$ Then 
$$\wide\,(P_2)=\add[(P_2/\rad P_2)\oplus (P_3/\rad^2 P_3)\oplus P_2\oplus P_3\oplus (P_2/\rad^3P_2)].$$
\end{ex}

In what follows, we give some discussion and a generalization of Theorem \ref{CR3}.

\begin{rk}\label{CR4} Let $M=\oplus_{i=1}^t M_i$ with $M_i$ indecomposable and $M_i\not\simeq M_j$ for $i\neq j.$
\begin{enumerate}
\item[(1)] If $\Ext_\Lambda^1(M,M)=0$ and $\End_\Lambda(M)$ is a product of division rings, then   $\add\,(M)=\wide\, (M).$

\noindent
{\rm{This follows directly from  Proposition \ref{R2}, since $\Ext_\Lambda^1(M,M)=0$ implies that $\add (M)$ is closed under extensions.}}
\item[(2)] The converse of {\rm (1)}  does not hold in general. 

\noindent
{\rm To see this, consider an artin algebra $\Lambda$  of finite representation type, which is not semi-simple and
such that the set of indecomposable $\Lambda$-modules (up to isomorphism)  is given by the set 
$\{M_1,M_2,\ldots, M_n\}$ with $n \geq 2.$ Take $M:=\oplus_{i=1}^nM_i.$ It is clear that $\add (M)=\modu (\Lambda)$ and therefore 
$\add (M)$ is wide, but $\End_\Lambda (M)$ is not  semi-simple neither $\Ext_\Lambda^1(M,M)= 0.$ }
\end{enumerate}
\end{rk}

The next result states that the converse of item (1)  in Remark \ref{CR4} holds whenever $M$ is a projective $\Lambda$-module.

\begin{pro}\label{CR5}
Let $P=\oplus_{i=1}^t P_i$ with each $P_i$ indecomposable projective $\Lambda$-module  such that $P_i\not\simeq P_j$ for $i\neq j.$ Then, the following statements are equivalent.
\begin{itemize}
\item[(a)]  $\add\,(P)=\wide\,(P).$
\item[(b)]   $\End_\Lambda (P)=\times_{i=1}^t\End_\Lambda(P_i)$ and every $\End_\Lambda (P_i)$ is a division ring.
\item[(c)]  $\End_\Lambda (P)$ is a product of division rings.
\end{itemize}
\end{pro}
\begin{dem}  We only show that the first and the second statements are equivalent. Clearly, the second statement implies the third. We leave for the reader to show that the third implies the second. 
\

Let us show the equivalence between the first the second statement.  We first assume  that  $\add\,(P)$ is a wide subcategory of $\modu\,(\Lambda).$
\

In order to get the decomposition $\End_\Lambda (P)=\times_{i=1}^t\End_\Lambda(P_i)$ as rings, it is enough to see that $\Hom_\Lambda(P_i,P_j)=0,$ for any $i\neq j.$ Indeed, suppose there exists a non-zero morphism $f:P_i\rightarrow P_j$ for some $i\neq j.$ Since $\add\,(P)$ is wide, we have that $0\to \Ker\,(f)\to P_i\to \Ima\,(f)\to 0$ and $0\to \Ima\,(f)\to P_j\to \Coker\,(f)\to 0$ are 
exact sequences in $\add\,(P),$ and thus both of them split. So, we get the decompositions 
$P_i=\Ker\,(f)\oplus \Ima\,(f)$ and $P_j=\Ima\,(f)\oplus \Coker\,(f).$ Using now that $P_i$ and $P_j$ 
are indecomposable and the fact that $\Ima\,(f)\neq 0,$ it follows that $P_i\simeq \Ima\,(f)=P_j;$ 
contradicting that $P_i\not\simeq P_j.$
\

We prove now that $\End_\Lambda (P_i)$ is a division ring. Indeed, let $g:P_i\to P_i$ be a non-zero morphism. Since $\add\,(P)$ is wide, we have that $0\to \Ker\,(g)\to P_i\to \Ima\,(g)\to 0$ is 
an exact sequence in $\add\,(P),$ and so it splits giving us $P_i=\Ker\,(g)\oplus \Ima\,(g).$ Since 
$P_i$ is indecomposable and $g\neq 0,$ it follows that $g$ is an isomorphism.

The other implication follows from item (1) of the former remark.
\end{dem}
 
 \section{Standard modules and wide subcategories}
 
Let $\Lambda$ be an artin $R$-algebra and $D:=\Hom_\Lambda(-,I):\modu\,(\Lambda)\to \modu\,(\Lambda^{op})$ be the usual duality functor, where $I$ is the injective envelope of $R/\rad\,(R).$
\

In this section we will  prove that the category $\wide\,(\Delta)$ coincides with $\modu (\Lambda)$, where $\Delta$ is the set of standard modules. Furthermore, we analyse the case when 
the category $\F(\Delta),$ of $\Delta$-filtered modules, is wide. In order to do that, we will need some preliminary results. Recall that $\topo\,(M)$ denotes the module $M/\rad\,(M),$  for any $M\in\modu\,(\Lambda).$ 

\begin{lem}\label{sumando}
Let $\X\subseteq\modu (\Lambda)$ be a wide subcategory and $M\in\X.$ If $\topo (M)$ is a direct summand of $\soc\,(M)$ then $\topo\,(M)\in\X.$ 
\end{lem}

\begin{dem} Let $\topo (M)$ be a direct summand of $\soc\,(M).$ Then there is a monomorphism 
$\mu:\topo\,(M)\to \soc\,(M).$ Consider now the natural projection $\pi:M\to \topo\,(M)$ and the 
inclusion $i:\soc\,(M)\to M.$ Therefore the morphism $f:=i\mu\pi:M\to M$ satisfies that 
$\Ima\,(f)=\topo\,(M);$ and since $\X$ is wide we get that $\topo\,(M)=\Ima\,(f)\in\X.$
\end{dem}

\begin{pro}\label{standard}
Let $\X$ be a wide subcategory of $\modu\,(\Lambda)$ and let 
$$\{M_1,M_2,\cdots, M_t\}\subseteq\X$$ 
be such that $\topo (M_i)$ is equal to the simple $S_i$ for all $i.$ If
each module $M_i$ has only composition factors among the simple $\Lambda$-modules $S_j$ with 
$j\leq i,$ then   $\{S_1,S_2,\cdots, S_t\}\subseteq\X.$
\end{pro}
\begin{dem}
We will proceed by induction on $t$. By Lemma \ref{sumando}, the result is clear for $t=1$.
Let us assume that $\{S_1,S_2,\cdots, S_{j-1}\}\subseteq\X.$ To prove that $S_j\in \X$ we will use induction on $\ell_\Lambda (M_j).$ If $\ell_\Lambda (M_j)=1$ then $S_j=M_j\in\X.$ 
\

Let $\ell_\Lambda (M_j)>1.$ Consider the exact sequence
$$\varepsilon:\quad 0\rightarrow \Tr_{\oplus_{i<j}P_i}(\soc\, M_j)\rightarrow M_j\rightarrow M_j'\rightarrow 0.$$

\noindent If $\Tr_{\oplus_{i<j}P_i}(\soc\, M_j)=0$ then $\soc (M_j)=S_j^m$ and by  Lemma \ref{sumando} we get $S_j\in\X.$
\

Assume that $\Tr_{\oplus_{i<j}P_i}(\soc\, M_j)\neq 0.$ Then  $\ell_\Lambda (M'_j)<\ell_\Lambda (M_j)$ with $M'_j$ having  composition factors among the simple $\Lambda$-modules  
$S_1,S_2,\cdots S_j.$ Furthermore,  since $\Tr_{\oplus_{i<j}P_i}(\soc\, M_j)\in\F(\{S_1,S_2,\cdots, S_{j-1}\})\subseteq \X,$  $M_j\in\X$ and $\X$ is wide, from the exact sequence $\varepsilon$ we get that $M'_j\in \X$. Thus, by induction, since $\ell_\Lambda (M'_j)<\ell_\Lambda (M_j),$ we conclude that $S_j\in\X.$ 
\end{dem}

\begin{teo}\label{delta}
For the set $\Delta$ of standard $\Lambda$-modules, we have that 
$$\wide\,(\Delta)=\modu\,(\Lambda).$$
\end{teo}
\begin{dem} As we know $\Delta=\{\Delta(1),\Delta(2),\cdots,\Delta(n)\}$ is the set of standard $\Lambda$-modules with respect to the natural order $\leq$ on the 
set $\{1,2,\cdots,n\},$ where $n$ is the number of simple $\Lambda$-modules (up to isomorphism). It is 
also well known, that $\topo\,(\Delta(i))=S_i$ and $\Delta(i)$ has composition factors amount the 
simple $\Lambda$-modules $S_j$ with $j\leq i.$ Then, by Proposition \ref{standard}, we get that 
$\wide\,(\Delta)$ contains all the simple $\Lambda$-modules; proving the result.
\end{dem}

\begin{cor}
Let $\F(\Delta)$ be an abelian subcategory of $\modu\,(\Lambda).$ Then $\Lambda$ is a quasi-hereditary  algebra and $\gld (\Lambda)\leq 1+\pd\, D(\Lambda_\Lambda),$ where $D$ is the usual duality functor. 
\end{cor}

\begin{dem} It is well known that $\F(\Delta)$ is closed under extensions, and hence, by 
Remark \ref{criterioA} we 
get that $\F(\Delta)$ is a wide subcategory of $\modu\,(\Lambda).$ Then, by Theorem \ref{delta} we conclude that $\F(\Delta)=\modu (\Lambda).$ Therefore  ${}_\Lambda\Lambda\in\F(\Delta)$ and $\F(\Delta)$ is closed under submodules. Thus, from   \cite[Proposition 3.21]{MMS3} we get that $\Lambda$ is quasi-hereditary  and $\gld (\Lambda)\leq 1+\pd\, T,$ where $T$ is the characteristic tilting $\Lambda$-module associated to $\Lambda.$ Furthermore \cite[Corollary 4]{R} gives us that 
$T=\oplus_{i=1}^n\,I_i,$ where $I_i$ is the injective envelope of the simple $S_i.$ In particular, 
we have that $\pd\,(T)=\pd\,D(\Lambda_\Lambda).$
\end{dem}

\vspace{.4cm}
 It is well known that  $\F(\Delta)$ is closed under extensions, direct summands and kernels of epimorphisms. A natural question is to give necessary and sufficient conditions to get that 
 $\F(\Delta)$ is a thick subcategory of $\modu\,(\Lambda).$
The following result characterizes, in the case of standardly stratified algebras, when 
the category $\F(\Delta)$ is closed under cokernels of monomorphisms. We denote by 
$\Pro^{<\infty}(\Lambda)$ the class of all 
$\Lambda$-modules of finite projective dimension.

\begin{pro} \label{FD=P}
Let $\Lambda$ be a standardly stratified algebra. Then,
$\F(\Delta)$ is a thick subcategory of $\modu\,(\Lambda)$ if and only if 
$\F(\Delta)=\Pro^{<\infty}(\Lambda).$    
\end{pro}
\begin{dem}
$(\Leftarrow)$ If $\F(\Delta)={\Pro}^{<\infty}(\Lambda)$ then the conclusion is trivial, since $\Pro^{<\infty}(\Lambda)$ is a thick subcategory of $\modu\,(\Lambda)$.
\
 
$(\Rightarrow)$ Assume that $\F(\Delta)$ is closed under cokernels of monomorphisms. In order to obtain that  $\F(\Delta)={\Pro}^{<\infty}(\Lambda),$  by using  \cite[Proposition 3.18]{MMS3}, it is enough to prove
that $(\add \,T)^{\wedge}\subseteq\F(\Delta),$ where $T$ is the characteristic tilting 
$\Lambda$-module. For the notation of $(\add \,T)^{\wedge},$ see \cite[p. 397]{MMS3}.
\

Let $M\in (\add \,T)^{\wedge}.$  Then we have an exact sequence
$$0\rightarrow T_k\to T_{k-1}\to\cdots \to T_0\rightarrow M\rightarrow 0$$
with $T_i\in \add \,T$ for $i=0,1,\cdots, k.$ Since $\add\,T\subseteq\F(\Delta)$ and $\F(\Delta)$ is closed under cokernels of monomorphisms, we get that $M\in\F(\Delta).$ Therefore $(\add \,T)^{\wedge}\subseteq\F(\Delta),$ proving the result.
\end{dem}

\vspace{.2cm}

It is well known that if $\Lambda$ is a quasi-hereditary algebra then its global dimension is finite. So we get the following corollary.

\begin{cor}
Let $\Lambda$ be a quasi-hereditary algebra. If $\F(\Delta)$ is a thick subcategory of $\modu\,(\Lambda)$
then $\F(\Delta)=\modu\,(\Lambda).$
\end{cor}
\begin{dem} It follows from Proposition \ref{FD=P}.
\end{dem}

\vspace{.4cm}

Next, we consider a quotient path $k$-algebra $\Lambda:=kQ/I$ and characterize, in this case,  
when the category $\F(\Delta)$ is wide. For doing so, we will  make use of the following well known definition.

\begin{defi}
Let $\Lambda=kQ/I$ be a quotient path $k$-algebra and $\leq$ be a linear order on the set of vertices 
$Q_0=\{1,2,\cdots, n\}.$ We say that the algebra $\Lambda$ is triangular, with respect to 
the partially ordered set $(Q_0,\leq),$ if $Q$ does not have oriented cycles and $\Hom_\Lambda(\Lambda e_i,\Lambda e_j)=0$ for $i<j.$ 
\end{defi}

We remark that the above definition is equivalent to say that the ordinary quiver of $\Lambda$ is {\bf directed}. The definition of directed quiver was introduced by S. Koenig in \cite{Ko}.
 
\begin{teo}\label{triangular}
Let $\Lambda=kQ/I$ be  a quotient path $k$-algebra, $\leq$ be a linear order on $Q_0$ and 
$\Delta$ be the set of standard $\Lambda$-modules (with respect to the linear order $\leq$ on $Q_0$). Then, the following statements are equivalent.
\begin{itemize}
\item[(a)] $\F(\Delta)=\modu (\Lambda).$
\item[(b)] $\F(\Delta)$ is wide.
\item[(c)] $\Lambda$ is triangular with respect to the partially ordered set $(Q_0,\leq).$
\item[(d)] $\F(\Delta)$ is an abelian subcategory of $\modu\,(\Lambda).$
\end{itemize} 
\end{teo}
\begin{dem} Let $Q_0=\{v_1, v_2,\cdots,v_n\}$ where $v_1<v_2<\cdots<v_n.$ For simplicity, we consider 
$\Delta(i):=\Delta(v_i),$ $S_i:=S_{v_i}$ and $P_i:=P_{v_i}.$
\

 $(a)\Leftrightarrow (b)$ It follows from the Theorem \ref{delta}. 
 \
 
 $(a)\Rightarrow(c)$  Since  $\F(\Delta)=\modu (\Lambda),$ we get that $\Delta(i)=S_i$ for any $i;$ otherwise $S_i$ could not be in $\F(\Delta)$ and so  $P_n=\Delta(n)=S_n.$ Therefore $v_n$ is a sink and then 
 $\Hom_\Lambda(\Lambda e_j,\Lambda e_n)\simeq e_j\Lambda e_n=0$ for any $j<n.$
We claim now that $v_{n-1}$ is a sink in $Q-\{v_n\}.$ Indeed, we know that $S_{n-1}=\Delta(n-1)=P_{n-1}/\Tr_{P_n}(P_{n-1}).$ Suppose that $v_{n-1}$ is not a sink in  $Q-\{v_n\}.$ Then there exists $t<n-1$ and an arrow $v_{n-1}\rightarrow v_t$ in $Q.$ This implies that the simple $S_{t}$ is a composition factor of $\Delta(n-1)=S_{n-1},$ which is a contradiction. Hence $v_{n-1}$ is a sink in  $Q-\{v_n\}.$
Inductively, we get that $v_{n-i}$ is a sink in  $Q-\{v_n,v_{n-1},\ldots, v_{n-i+1}\}.$ Therefore $Q$ does not have oriented cycles and $\Hom_\Lambda(\Lambda e_i,\Lambda e_j)\simeq e_i\Lambda e_j=0$ for all $i<j.$ 
\

$(c)\Rightarrow(a)$ Since $Q$ does not have oriented cycles and  $\Hom_\Lambda(\Lambda e_i,\Lambda e_j)=0$ for $i<j,$ we get that 
$\Tr_{\oplus_{j>i} P_j}(P_i)= \Tr_{\oplus_{j\neq i} P_j}(P_i).$ Therefore $$\Delta(i)=P_i/\Tr_{\oplus_{j>i} P_j}(P_i)=P_i/\Tr_{\oplus_{j\neq i} P_j}(P_i)\simeq S_i.$$
Then, we get that 
$$\F(\Delta(1),\Delta(2),\cdots,\Delta(n))=\F(S_1,S_2,\cdots,S_n)=\modu (\Lambda).$$
Finally, from Theorem \ref{delta}  and the  using the fact that $\F(\Delta)$ is closed under extensions, we get that the items (b) and (d) are equivalent.
\end{dem}
\vspace{0.2cm}

Let $\Lambda$ be an artin $R$-algebra. We finish this section with the following proposition.

\begin{pro}
Let $\X$ be a wide subcategory of $\modu (\Lambda)$ and $I$ a two sided ideal of $\Lambda$   such that $\Lambda/I\in\X.$ Then $\soc (\Lambda/I)\in\X.$ 
\end{pro}

\begin{dem} Let us see first that 
$$ (*) \quad \forall\,  t\in\Lambda /I\quad\ann_{\Lambda/I}(t)\in\X.$$
Indeed, for any $t\in\Lambda,$ the morphism $\phi_t:\Lambda/I\rightarrow \Lambda/I$, given by $\phi_t(x):=xt$ satisfies $\ann_{\Lambda/I}(t)=\Ker (\phi_t)\in\X,$ proving $(*).$
\

Since $\soc (\Lambda /I)$ is a finitely generated $R$-module, we have that  $\soc (\Lambda /I)=\langle x_1,x_2,\ldots, x_m\rangle_R.$ Thus, from  $(*)$ and Proposition \ref{R1}, we get that 
$$\ann(\soc(\Lambda /I))=\cap_{i=1}^m \ann_{\Lambda/I}(x_i)\in \X$$ 
and therefore $(\Lambda/I)/(\ann(\soc(\Lambda/I)))\in\X.$
Finally, since 
$$\add\, (\soc(\Lambda /I))=\add\, [(\Lambda/I)/\ann (\soc(\Lambda/I))],$$
 we get that $\soc (\Lambda/I)\in\X.$
\end{dem}

\section{Wide categories of type $\pres\,(P)$ with $P$ projective}

Let $\Lambda$ be an artin $R$-algebra. In this section, we characterise when a wide subcategory 
$\X$ of $\modu\,(\Lambda)$ is of the form $\pres\,(P)$ for some $P\in\proj\,(\Lambda).$
\

Let $P\in\proj\,(\Lambda).$ We recall, for details see  \cite{ARS}, that $\pres\,(P)$ is the full subcategory of 
$\modu\,(\Lambda),$ whose objects are all the $\Lambda$-modules $M$ admitting  a presentation in $\add\,(P),$ that is,  an exact sequence 
$P_1\to P_0\to M\to 0$ with $P_0,P_1\in\add\,(P).$ Using the  Horseshoe Lemma, we see  that $\pres\,(P)$ is closed under extensions.
\

Let $\Gamma:=\End_\Lambda(P)^{op}.$ In this case, it is well known, see \cite[Proposition 2.5]{ARS}, that the evaluation functor $(P,-):\pres\,(P)\to \modu\,(\Gamma)$ is an equivalence of categories. So, 
by using this equivalence, we can translate the abelian structure from $\modu\,(\Gamma)$ to $\pres\,(P).$ In general, it can happen that $\pres\,(P)$ is not an abelian subcategory of $\modu\,(\Lambda),$ 
although $\pres\,(P)$ is an abelian category. In what follows, we give an example illustrating this 
situation.

\begin{ex}
Let    $\Lambda=k\tilde{A}_2/I$  where $\tilde{A}_2$ is the quiver

$$
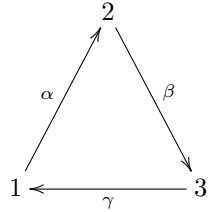
\xymatrix{
						&2 \ar@{->}[ddr]^\beta 	&\\
						&						&\\
1 \ar@{->}[uur]^\alpha 	&						&3 \ar@{->}[ll]^\gamma
}
$$
and $I=J^5.$ Recall that $J$ denotes de arrow ideal in $k\tilde{A}_2.$ Let $P:=P_1$ be the projective $\Lambda$-module attached to the vertex $1.$ We assert that $\pres\,(P)$ is not an abelian subcategory of $\modu\,(\Lambda).$ Observe that in this case 
$\pres\,(P)$ is
equivalent to $\modu\, (k[x]/(x^2)).$
\

{\rm Indeed, consider the natural projection $\pi:P\to P/\rad^2(P).$ Note that 
$P/\rad^2(P)\simeq \rad^3(P)$ and hence we have a monomorphism $\mu: P/\rad^2(P)\to P.$ Let 
$g:=\mu\pi:P\to P.$ Note that $\Ker\,(g)\simeq P_3/\rad^3(P_3)\not\in\add\,(P).$ Therefore 
$\Ima\,(g)\not\in\pres\,(P);$ proving that $\pres(P)$ is not an abelian subcategory of $\modu\,(\Lambda).$
}
\end{ex}

For a given subcategory $\X\subseteq\modu\,(\Lambda),$ we consider the class
 $\mathbb{P}_0(\X):=\{P_0(X)\;:\:X\in\X\},$ where $P_0(X)$ is the projective cover of $X.$ We recall that, an exact sequence in $\X$ is just an exact sequence in $\modu\,(\Lambda),$ whose terms are all in $\X.$ 
 \
 
 Recall that $\proj\,(\X)$ denotes the full subcategory of $\X$ whose objects are relatively projective, that is, $P\in\proj\,(\X)$ if and only if for any exact sequence $0\to A\to B\stackrel{f}{\to} C\to 0$ in $\X,$ the 
 map $\Hom_\Lambda(P,f):\Hom_\Lambda(P,B)\to\Hom_\Lambda(P,C)$ is surjective. We say that $\X$ has enough projectives if any $M\in\X$ admits an exact sequence $0\to K\to P\to M\to 0$ in $\X$ with 
 $P\in\proj\,(\X).$

\begin{teo}\label{wideProj} Let $\X$ be a wide subcategory of $\modu\,(\Lambda).$ Then, the following statements are equivalent.
\begin{itemize}
\item[(a)] $\X=\pres\,(P)$ for some $P\in\proj\,(\Lambda).$
\item[(b)] $\X$ has enough projectives and $\proj\,(\X)\subseteq \proj\,(\Lambda).$
\item[(c)] $\X$ is closed under projective covers, that is $\mathbb{P}_0(\X)\subseteq\X.$
\end{itemize}
If one of the above equivalent conditions holds, then 
$$\add\,(P)=\proj\,(\X)=\add\,(\mathbb{P}_0(\X)).$$ 
In particular, $P$ is uniquely determined  up to additive closures, and  the functor $\Hom_\Lambda(P,-):\X\to 
\modu\,(\End_\Lambda(P)^{op})$ is an equivalence of categories.
\end{teo}
\begin{dem}
$(a)\Rightarrow(b)$ It is clear that $\add\,(P)\subseteq \proj\,(\X).$ We assert that, for any 
 $M\in\pres\,(P),$ there is an exact sequence $0\to K\to P_0\to M\to 0$ in $\X,$ with 
$P_0\in\add\,(P).$ 
\

Indeed, let $M\in\pres\,(P).$ Then there is an exact sequence $P_1\stackrel{f}{\to}P_0\to M\to 0,$ 
with $P_0,P_1\in\add\,(P).$ Using the fact that $\X$ is wide, we get that $\Ima\,(f)\in\pres\,(P);$ and hence $0\to \Ima\,(f)\to P_0\to M\to 0$ is the desired exact sequence. This, in particular, 
implies that $\X$ has enough projectives, since $\add\,(P)\subseteq \proj\,(\X).$
\

Let us prove that $\proj\,(\X)\subseteq \add\,(P).$ Consider $Q\in \proj\,(\X).$ Then, there is an 
exact sequence $\eta:\quad 0\to K\to P_0\to Q\to 0$ in $\X,$ with 
$P_0\in\add\,(P).$ Thus, $\eta$ splits and then $Q\in\add\,(P).$
\

$(b)\Rightarrow(a)$ Let $P:=\oplus_{i=1}^tP_i,$ where  $\{P_1,P_2,\cdots,P_t\}$ is a 
complete list of indecomposable pairwise non isomorphic $\Lambda$-modules in $\proj\,(\X)\subseteq \proj\,(\Lambda).$ Furthermore, $\add\,(P)\subseteq \X,$ since $P\in\X.$
Let $M\in\pres\,(P).$ Then there is an exact sequence $P_1\stackrel{f}{\to}P_0\to M\to 0,$ 
with $P_0,P_1\in\add\,(P).$ Since $\add\,(P)\subseteq \X,$ it follows that $M\simeq \Coker\,(f)\in\X;$ proving that $\pres\,(P)\subseteq\X.$
\

Let $M\in\X.$ Using that $\X$ has enough projectives and $\proj\,(\X)=\add\,(P),$ we get an 
exact sequence $0\to K\to Q_0\to M\to 0$ in $\X,$ with $Q_0\in\add\,(P).$ Doing the same with $K,$ 
we get an exact sequence $Q_1\to Q_0\to M\to 0,$ with $Q_0,Q_1\in\add\,(P).$ Therefore 
$M\in\pres\,(P)$ and thus $\X\subseteq \pres\,(P).$
\

$(b)\Rightarrow(c)$ Let $M\in\X.$ Then there is an epimorphism $Q\to M,$ with $Q\in\proj\,(\X)\subseteq\proj\,(\Lambda).$ Therefore, the projective cover $P_0(M)$ is a direct summand of $Q$ and 
thus $P_0(M)\in\X.$
\

$(c)\Rightarrow(b)$ Since $\mathbb{P}_0(\X)\subseteq\X$ and $\X$ is closed under kernels of epimorphism, we have an exact sequence $0\to K_M\to P_0(M)\to M\to 0$ in $\X,$ for any $M\in\X.$ By 
taking $M:=Q\in\proj\,(\X),$ we get that $Q$ is a direct summand of $P_0(M)$ and thus $Q\in\proj\,(\Lambda).$ Finally, it is clear that $P_0(M)\in\proj\,(\X)$ and then $\X$ has enough projectives.
\end{dem}

Let $\varphi:A\to B$ be a ring epimorphism of finite dimensional $k$-algebras, and let $F_\varphi:\modu\,(B)\to\modu\,(A)$ be 
its  associated restriction functor. It is well known that $\modu\,(B)$ can be identified, as a full subcategory of $\modu\,(A),$  with 
the essential image of the functor $F_\varphi.$ Let $\X$ be a class  of finitely generated $A$-modules. We say  that $\X$ is precovering  (or contravariantly finite) if for any $M\in\modu\,(\Lambda),$ there is a morphisms $f:X\to M,$  with $X\in\X,$ such that $(Z,f):(Z,X)\to (Z,M)$ is bijective for any 
$Z\in\X.$ Dually, there is the notion of preenveloping (or covariantly finite) class. The class $\X$ is functorially finite if it is precovering and preenveloping. 
\

The following proposition appears in \cite[Proposition 4.1]{MS} and it is essentially 
done in  \cite{GL} (see also in \cite[Theorem 1.6.1]{I}). We 
 recall that two ring epimorphisms $f : A \to B$ and $g : A \to C$ are equivalent if there is a (necessarily unique) isomorphism of rings
$h : B\to  C$  such that $g = h f .$

\begin{pro}\label{ringep1}  Let $A$ be a finite dimensional $k$-algebra. By assigning to a given ring epimorphism,
the essential image of its associated restriction functor, one gets a bijection between
\begin{itemize}
\item[(a)] equivalence classes of ring epimorphisms $A\to B$ with $\dim_k(B) $ finite and $\mathrm{Tor}^{A}_1(B,B)=0;$
\item[(b)]  functorially finite wide subcategories of $\modu\,(A).$
\end{itemize}
\end{pro}

In case of wide functorially finite subcategories, we have that Theorem \ref{wideProj} can be connected, with the proposition above, as can be seen in the following corollary.

\begin{cor}\label{ringepi2} Let $A$ be a finite dimensional $k$-algebra. Then, for  a functorially finite wide subcategory $\X$ of $\modu\,(A)$  the following statements are equivalent.
\begin{itemize}
\item[(a)] There is a ring epimorphism $A\to B$ such that $B\in\proj\,(A)$ and $\X=\pres\,({}_AB).$
\item[(b)] $\X$ has enough projectives and $\proj\,(\X)\subseteq \proj\,(A).$
\end{itemize}
If one of the above equivalent conditions hold,  there is a basic  $A$-module $P$ such that $\add\,(P)=\proj\,(\X)$ and $B$ is 
Morita equivalent to $\End_A(P)^{op}.$
\end{cor}
\begin{dem} (a) $\Rightarrow$ (b)  Suppose (a) holds. Then by Proposition \ref{ringep1}, we get $\X=\modu\,(B)$ and so $\X$ has enough projectives. Furthermore, any indecomposable projective $B$-module is a direct summand of ${}_BB,$ and hence projective as an $A$-module as $B$ is so. 
\

(b) $\Rightarrow$ (a) Assume that (b) holds. by Proposition \ref{ringep1} there is a ring epimorphism $A\to B$ with $\dim_k(B)$ finite 
and $\X=\modu\,(B).$ But then $B\in\proj\,(\X)\subseteq \proj\,(A).$
\

Finally, the last assertion follows from Theorem \ref{wideProj}.

\end{dem}

\footnotesize

\vskip3mm \noindent Eduardo Marcos:\\
Instituto de Matem\'aticas y Estadistica,\\
Universidad de Sao Paulo,\\
Sao Paulo, BRASIL.

{\tt enmarcos@gmail.com}

\vskip3mm \noindent Octavio Mendoza:\\
Instituto de Matem\'aticas,\\
Universidad Nacional Aut\'onoma de M\'exico,\\
Circuito Exterior, Ciudad Universitaria,\\
M\'exico D.F. 04510, M\'EXICO.

{\tt omendoza@matem.unam.mx}

\vskip3mm \noindent Corina S\'aenz:\\
Departamento de Matem\'aticas, Facultad de Ciencias,\\
Universidad Nacional Aut\'onoma de M\'exico,\\
Circuito Exterior, Ciudad Universitaria,\\
M\'exico D.F. 04510, M\'EXICO.

{\tt corina.saenz@gmail.com}

\vskip3mm \noindent Valente Santiago:\\
Departamento de Matem\'aticas, Facultad de Ciencias,\\
Universidad Nacional Aut\'onoma de M\'exico,\\
Circuito Exterior, Ciudad Universitaria,\\
M\'exico D.F. 04510, M\'EXICO.

{\tt valente.santiago.v@gmail.com}


\vspace{1cm}





\begin{thebibliography}{60}
\bibitem{AHLU} \'Agoston, I.; Happel, D.; Luk\'acs, E.; and Unger, L.  Standardly stratified algebras and tilting. {\it {J. Algebra.}} 226, 
144-160, (2000).
\bibitem{ARS} Auslander, M.; Reiten, I.; Smal\o, S.O. (1995). Representation theory of Artin algebras. {\it {Cambridge University Press}}.

\bibitem{ASS}  Assem, I.; Simson, D.; Skowro\'nski, A.(2006). Elements of the Representation Theory of Associative Algebras. Volumen 1: Techniques of Representation Theory. {\it{LMS, Student text}} 65, Cambridge.

\bibitem{CPS} Cline, E.; Parshall, B.J.; Scott, L.L. Finite dimensional algebras and highest weight categories. {\it {J. Reine Angew. Math.}} 391,  85-99, (1988).

\bibitem{GL} Geigle, W.; Lenzing, H. Perpendicular categories with applications to representations and sheaves. {\it{ J. Algebra.}}
144, no. 2, 273-343, (1991).

\bibitem{H} Hovey, M. Classifying subcategories of modules. {\it{Trans. of the AMS}} 353 (8), 3181-3191, (2001).

\bibitem{IT} Ingalls, C.; Thomas, H. Noncrossing partitions and representations of quivers. {\it{Compositio Math.}} 145, 1533-1562, (2009).
\bibitem{IPT} Ingalls, C.; Paquette, C.; Thomas H. Semi-stable subcategories for Euclidean quivers. arXiv:1212.1424v3.

\bibitem{I} Iyama, O. Rejective subcategories of artin algebras and orders, preprint, arXiv:math/0311281.

\bibitem{Ko} Koenig, S. Exact Borel subalgebras of quasi-hereditary algebras I. With an appendix by Leonard Scott. {\it{Math. Z.}} 220, no. 3, 399-426, (1995).

\bibitem{K} Krause, H. Thick subcategories of modules over commutative noetherian rings (with an appendix by Srikanth Iyengar). {\it{Math. Ann.}} 340, 733-747, (2008).

\bibitem{MMS3} Marcos, E.N.; Mendoza, O.; S\'aenz, C. Applications of Stratifying Systems to the finitistic dimension. {\it{Journal of pure and 
applied algebra}}  205 (2), 393-411, (2006).

\bibitem{MS} Marks F.; Stovicek J. Torsion classes, wide subcategories and localisations. arXiv:1503.04639v1, (2015).

\bibitem{P} Pogorzaly, Z. On bimodules determining stable equivalences for self-injective algebras. Lecture given in ARTA IV, CIMAT, Guanajuato, M\'exico, (2015).


\bibitem{R} Ringel, C.M. The category of modules with good filtrations over a quasi-hereditary algebra has almost split sequences.
{\it{Math. Z.}} 208, 209-233, (1991).

\bibitem{S} Scott, L.L. Simulating algebraic geometry with algebra. I. The algebraic theory of derived categories, in 
{\it{The Arcata Conference on Representations of Finite Groups (Arcata, Calif., 1986)}}, 271-281, Proc. Sympos. Pure Math., 47, Part 2, Amer. Math. Soc., Providence, RI, (1987).

\bibitem{T} Takahashi, R. Classifying subcategories of modules over a commutative noetherian ring. 
{\it{Journal of the LMS}} 78 (3), 767-782, (2008).

\end{thebibliography}
\end{document}